# A Data-Driven Warm Start Approach for Convex Relaxation in Optimal Gas Flow

Haizhou Liu, Lun Yang, Xinwei Shen, *Member, IEEE*, Qinglai Guo, *Senior Member, IEEE*, Hongbin Sun, *Fellow, IEEE* and Mohammad Shahidehpour, *Fellow, IEEE*

*Abstract*— In this letter, we propose a data-driven warm start approach, empowered by artificial neural networks, to boost the efficiency of convex relaxations in optimal gas flow. Case studies show that this approach significantly decreases the number of iterations for the convex-concave procedure algorithm, and optimality and feasibility of the solution can still be guaranteed.

## Nomenclature

*A. Indices and Sets*

| | |
|---|---|
| $s \in \Omega_S$ | Set of sources. |
| $g \in \Omega_N$ | Set of nodes. |
| $(m,n) \in \Omega_P$ | Set of pipelines, with $m/n$ denoting the start/end node. |
| $(i,j) \in \Omega_C$ | Set of compressors, with $i/j$ denoting the start/end node. |
| $t \in \Omega_T$ | Set of time slots. |

*B. Variables*

| | |
|---|---|
| $G_s$ | Gas output of source $s$. |
| $F_{mn}, F_{m'n'}$ | Gas flow of pipelines $(m,n)$, $(m',n')$. |
| $\pi_g, \pi_i, \pi_j, \pi_m, \pi_n$ | Nodal pressure of gas nodes $g$, $i$, $j$, $m$ and $n$. |
| $F_{C,ij}, F_{C,i'j'}$ | Gas flow in compressors $(i,j)$, $(i',j')$. |
| $W_{ij}$ | Gas consumption of compressor $(i,j)$. |
| $q_{mn}^{in}, q_{mn}^{out}$ | Gas flow into and out of pipeline $(m,n)$. |
| $M_{mn}$ | Natural gas stored in pipeline $(m,n)$. |
| $M_{mn,0}, M_{mn,T}$ | Natural gas stored in pipeline $(m,n)$ at the beginning and end of the time range $\Omega_T$. |
| $(\,)_t, (\,)_{t-1}$ | Values of variable $(\,)$ at time slots $t$ and $t-1$. |

*C. Parameters*

| | |
|---|---|
| $C_s$ | The unit gas production cost of source index $s$. |
| $\gamma_{ij}$ | The gas consumption ratio of compressor $(i,j)$. |
| $R_{ij}$ | Maximum compressibility ratio of compressor $(i,j)$. |
| $C_{mn,0}$ | Transmission parameter of pipeline $(m,n)$. |
| $F_{mn,0}$ | Warm Start of the gas flow of pipeline $(m,n)$. |
| $\lambda_g$ | Nodal load fluctuations at node $g$. |
| $L_g$ | The base nodal gas load at node $g$. |
| $F_{mn}^{max}$ | Maximum gas flow in pipeline $(m,n)$. |
| $G_s^{min}, G_s^{max}$ | Minimum and maximum gas output of source $s$. |
| $\pi_g^{min}, \pi_g^{max}$ | Minimum and maximum pressure of node $g$. |
| $F_{C,ij}^{max}$ | Maximum gas consumption of compressor $(i,j)$. |

## I. Introduction

THE increasing number of gas-fired generators and power-to-gas units has significantly intensified the coupling of electricity and natural gas systems. However, the gas flow equation of the natural gas system is highly nonconvex, which in turn complicates the optimization of the electricity-gas coupled system. Many algorithms have emerged to address the gas flow equation. For example, piecewise linearization is introduced to approximate gas flow [1]. However, many binaries are introduced, which would increase the computational burden. Nonlinear methods such as primal-dual interior point are also adopted, but their convergence cannot be guaranteed. Recently, convex relaxation, featuring second-order cone relaxation (SOCR) has shown superiority. SOCR includes many variants such as convex-concave procedure (CCP) [2] and convex envelope [3]. Many SOCRs, however, are either inaccurate in relaxation, or require many iterations to reach optimality. Additional improvements to SOCR-based optimization are therefore necessary.

On the other hand, artificial neural network (ANN) has become a powerful tool in supervised learning. ANN learns the nonlinear input-output mapping from massive data streams, to predict highly efficient and accurate outputs. Inspired by its distinct properties, we propose a data-driven convex relaxation approach for the optimal gas flow problem. ANN provides an accurate prediction of natural gas pressures, which is set as the warm start for SOCR techniques such as CCP. Case studies show that the proposed method can significantly decrease the required number of iterations for convergence, while guaranteeing the optimal solution for optimal gas flow.

## II. Model Formulation

### A. Steady-state Optimal Gas Flow

For steady-state optimal gas flow, the optimization model can be completely decoupled in each time slot. The objective function is to minimize the total gas output cost (1). (2)-(6) are upper/lower bound constraints for nodal gas pressure, gas flow, transported gas flow through compressor, gas source outputs and compressor ratio. The gas consumption of compressors is modeled by (7) [3]. The Weymouth gas flow equation is applied in (8). The nodal balance equation is given in (9).

$$\text{minimize} \sum_{s \in \Omega_S} C_s G_s \quad (1)$$

$$\pi_g^{min} \leq \pi_g \leq \pi_g^{max}, \forall g \in \Omega_N \quad (2)$$

$$-F_{mn}^{max} \leq F_{mn} \leq F_{mn}^{max}, \forall (m,n) \in \Omega_P \quad (3)$$

This work is supported by the National Key Research and Development Program SQ 2020YFE0200400, National Natural Science Foundation of China (No. 52007123) and the Science, Technology and Innovation Commission of Shenzhen Municipality (No. JCYJ20170411152331932).

Haizhou Liu, Lun Yang and Xinwei Shen are with Tsinghua-Berkeley Shenzhen Institute, Tsinghua Shenzhen International Graduate School, Tsinghua University. Qinglai Guo and Hongbin Sun are with Dept. of Electrical Engineering, Tsinghua University. Mohammad Shahidehpour is with the Electrical and Computer Engineering Department, Illinois Institute of Technology. (Corresponding author: Xinwei Shen, email: sxw.tbsi@sz.tsinghua.edu.cn; Hongbin Sun, email: shb@tsinghua.edu.cn).

$$0 \le F_{C,ij} \le F_{C,ij}^{\max}, \forall (i,j) \in \Omega_C \tag{4}$$

$$G_s^{\min} \le G_s \le G_s^{\max}, \forall s \in \Omega_S \tag{5}$$

$$\pi_i \le \pi_j \le R_{ij} \cdot \pi_i, \forall (i,j) \in \Omega_C \tag{6}$$

$$\tau_{ij} = \gamma_{ij} F_{C,ij}, \forall (i,j) \in \Omega_C \tag{7}$$

$$F_{mn}|F_{mn}| = C_{mn}^2 (\pi_m^2 - \pi_n^2), \forall (m,n) \in \Omega_P \tag{8}$$

$$\sum_{s \in g} G_s - \sum_{\substack{m \in g \\ (m,n) \in \Omega_P}} F_{mn} + \sum_{\substack{n' \in g \\ (m',n') \in \Omega_P}} F_{m'n'} - \sum_{\substack{i \in g \\ (i,j) \in \Omega_C}} (F_{C,ij} + \tau_{ij}) + \sum_{\substack{j' \in g \\ (i',j') \in \Omega_C}} F_{C,i'j'} = L_g, \forall g \in \Omega_N \tag{9}$$

### B. Quasi-Dynamic Optimal Gas Flow

In the quasi-dynamic optimal gas flow problem, the stored linepack is assumed to be adjustable with time, which requires a joint optimization of multiple time periods. The new objective is the total cost across the time range $\Omega_T = \{1,2,\ldots,T\}$

$$\text{minimize} \sum_{s \in \Omega_S, t \in \Omega_T} C_s G_{s,t} \tag{10}$$

Inequalities (2)-(6) and equalities (7)-(8) still hold, which are applied to all the time slots by adding a time $t$ as subscript to all variables. (9) is also modified into (11) to consider different gas inflow and outflow:

$$\sum_{s \in g} G_{s,t} - \sum_{m \in g, (m,n) \in \Omega_P} q_{mn}^{in} + \sum_{n' \in g, (m',n') \in \Omega_P} q_{m'n'}^{out} - \sum_{i \in g, (i,j) \in \Omega_C} (F_{C,ij,t} + W_{ij,t}) + \sum_{j' \in g, (i',j') \in \Omega_C} CF_{i'j',t} = \lambda_{g,t} L_{g,t}, \forall g \in \Omega_N, t \in \Omega_T \tag{11}$$

In addition, (12) defines the gas flow as the average of gas inflows and outflows. (13) relates linepack to nodal gas pressures. (14) specifies linepack changes with time. (15) states that total linepack in the initial and final states must be equal.

$$F_{mn,t} = \frac{1}{2}(q_{mn,t}^{out} + q_{mn,t}^{in}), \forall (m,n) \in \Omega_P, t \in \Omega_T \tag{12}$$

$$M_{mn,t} = H_{mn} \cdot \frac{1}{2}(\pi_{m,t} + \pi_{n,t}), \forall (m,n) \in \Omega_P, t \in \Omega_T \tag{13}$$

$$M_{mn,t} = M_{mn,t-1} + q_{mn,t}^{in} - q_{mn,t}^{out}, \forall (m,n) \in \Omega_P, t \in \Omega_T \tag{14}$$

$$\sum_{(m,n) \in \Omega_P} M_{mn,0} = \sum_{(m,n) \in \Omega_P} M_{mn,T} \tag{15}$$

## III. DATA-DRIVEN WARM STARTS

Despite the impressive performance of SOCR-based convex relaxation, including CCP as confirmed in [2] and [4], there are two unresolved issues: (1) A set of accurate initial points would need to be provided, which is crucial to the convergence of SOCR methods [2]; (2) The gas flow direction would need to be determined. The binary variable introduced for describing unknown gas flow directions will limit the computation efficiency of the SOCR methods.

Data-driven methods, empowered by ANN, can overcome these shortcomings by providing a proper warm start. To illustrate, a two-stage data-driven warm start method for CCP is given in Algorithm 1. In stage 1, an ANN is constructed using a multi-layer feedforward architecture, which is stated as

$$x_k^n = \text{ReLU}\left(\sum_{i=1}^{N_{n-1}} \omega_{i,k}^{n-1} x_i^{n-1} + b_k^n\right) \tag{16}$$

where $\text{ReLU}(\cdot)$ is the activation function and $x_k^n$ stands for the $k$-the neuron at the $n$-th layer. The weights $\omega_{i,k}^{n-1}$ and biases $b_k^n$ of the neural network are initialized randomly and updated via gradient descent in the training phase, using the training set which is obtained from a model-driven optimization method. The mean-squared error is chosen as the metric to evaluate the ANN's performance. To suppress gradient swaying, RMSprop is applied to consider the moving average of the gradients, e.g.

$$\omega := \omega - \eta \varphi / \sqrt{v + \varepsilon} \tag{17}$$

where $\eta$ and $\varepsilon$ are parameters. $\varphi$ denotes the magnitude of the gradient, and $v$ represents its moving average.

In stage 2, the as-trained ANN is used to initialize the pressures $\pi_m, \pi_n$ of new load scenarios. The direction of gas flow can also be inferred by the comparison of nodal pressures. Then CCP is iteratively carried out, with the initialized nodal pressures as warm starts. A non-negative slack variable $s_{mn}$, which is penalized, is introduced to allow mild violations of the gas flow equation at early stages of iteration. The iterative process is terminated once the gas flow violation is less than a threshold $\xi_0$. Algorithm 1 only shows the steady-state case, but can be converted to the quasi-dynamic case by a simple analogy.

---

**Algorithm 1** Data-Driven Warm-Start for CCP (Steady-State Optimal Gas Flow)

**Stage 1: Construct and train ANN.**
1. Construct an ANN with feedforward recursion (16). Randomly initialize $\omega$ and $b$.
2. Obtain the training set for ANN by pre-solving massive model-driven scenarios.
3. Train the ANN using
   (a) massive training set scenarios as inputs,
   (b) RMSProp (17) as back-propagation mechanism, and
   (c) mean squared error as metric.

**Stage 2: Find a warm start using the as-trained ANN, with which the CCP is iteratively carried out.**

4. For a new load scenario,
5. Initialize $\pi_m, \pi_n$ for each pipeline $(m,n)$ using the ANN. Denote the initialization values as $\pi_{m0}, \pi_{n0}$. Calculate $F_{mn,0}$ and infer gas flow direction using (8). Initialize $\tau = \tau_1$.
6. Optimize the following model (Suppose $\pi_m > \pi_n$):

   minimize $\sum_{s \in \Omega_S} C_s G_s + \sum_{(m,n) \in \Omega_P} \tau s_{mn}$

   subject to (2)-(9), $\left\|\begin{array}{c}\pi_n \\ F_{mn}/C_{mn}\end{array}\right\|_2 \le \pi_m, \forall (m,n) \in \Omega_P$, $s_{mn} \ge 0$ and

   $\pi_m^2 - (2\pi_{n0}\pi_n - \pi_{n0}^2) - \frac{1}{C_{mn}^2}(F_{mn,0}^2 - 2F_{mn}F_{mn,0}) \le s_{mn}, \forall (m,n) \in \Omega_P$.

7. Obtain the maximum violation of the gas flow constraint
   $\xi = \max_{(m,n) \in \Omega_P}\left\{\left|\sqrt{|\pi_m^2 - \pi_n^2|}/|F_{mn}| - 1\right|\right\}$.
8. If $\xi < \xi_0$
   Export the optimal solution. Finish.
   Else
   Update $\pi_{m0} := \pi_m, \pi_{n0} := \pi_n$. Calculate $F_{mn,0}$.
   $\tau = \min(\kappa\tau_1, \tau_{\max})$. Return to Step 6.

---

## IV. CASE STUDY

The data-driven CCP algorithm is implemented to analyze steady-state and quasi-dynamic solutions of the optimal gas flow of a 7-node and a 20-node natural gas system. Network parameters are provided in [5]. Table I lists the parameter settings for the data-driven CCP. In all scenarios, each nodal load is randomly fluctuating within ±10% around its base value.

TABLE I
PARAMETER SETTINGS FOR THE DATA-DRIVEN CCP.

| $\xi_0$ | $\rho$ | $\tau_1$ | $\tau_{\max}$ | $\kappa$ | $\varepsilon$ | $\eta$ |
|---|---|---|---|---|---|---|
| $10^{-3}$ | 1 | 1 | 1000 | 2 | $1 \times 10^{-8}$ | $10^{-2}$ |

For the physical model, the pressure range in the training set reaches an average of 26.1 Psig (pounds per square inch gauge) in the 7-node system. This suggests that providing a simple statistic (e.g., mean) of the pressure range as warm start might not be sufficient. However, ANN provides an accurate prediction of the nodal pressure $p_n$, given the nodal load fluctuations as input. The mean absolute error (MAE) of the test set reaches an average of 1.04 Psig, only 22% that of a dummy mean predictor (which always predicts the mean value of the pressure's upper and lower bounds).

To illustrate the effectiveness of the data-driven CCP, the following methods are compared.

M1: 100-segment piecewise linearization.

M2: Mixed-integer nonlinear solver SCIP [6].

M3: Cold-start CCP. Initial points are randomly chosen within the physical bounds. (Here a pre-determined gas flow direction is assumed to address (8)).

M4: The proposed data-driven CCP as in Algorithm 1.

All computation is conducted on a desktop computer with an Intel Core i5-7500 CPU and an 8G of RAM. The optimal results generated by the four methods are summarized in Table II. The data-driven CCP (M4) uses an ANN-based warm start, to converge the solution with a negligible optimality gap as compared to that of M1. Meanwhile, violation of the gas flow constraint $\xi$ is controlled down to $7.61 \times 10^{-6}$. On the other hand, for the cold-start CCP (M3), although a similarly low $\xi$ is reached, the obtained solutions are suboptimal, shown from the optimal cost with an optimality gap of 1.1%.

The computational efficiency of these methods is compared on both the 7-node and 20-node gas system, where the average performance under 50 random scenarios is summarized in Table III. In the 7-node gas system, thanks to a data-driven warm start, M4 requires an average of only 1.06 iterations to finish optimization, making it the most efficient among the 4 methods. Superiority in computational efficiency is more distinct in the 20-node gas system, where the average number of CCP iterations needed has significantly decreased from 14.17 down to 1.47 in the presence of warm starts, leading to a shorter computation time which is only 7.0% that of a cold-start CCP (M3). M4 also outperforms M1 and M2 by 10.7 and 3.9 times, respectively.

For the quasi-dynamic optimal gas flow, a 6-hour optimal gas flow problem is constructed on the 20-node natural gas system. An individual ANN is trained for each natural gas source at each time slot. The ANN-based data driven prediction method reaches an average prediction MAE of 1.72 Psig, which evaluates to a mere 1.45% of the average pressure.

The average accuracy and efficiency performances are listed in Table IV. The multiple time periods have rendered piecewise linearization (M1) inferior to other algorithms, both in terms of accuracy (indicated by a large $\xi$) and computational efficiency. This is also why we set M2 as benchmark for calculating the

TABLE II
DIFFERENT OPTIMIZATION RESULTS FOR THE STEADY-STATE OPTIMAL GAS FLOW OF A 7-NODE NATURAL GAS SYSTEM.

|    | Optimal Cost ($) | Source 1 Output (kcf*) | Source 2 Output (kcf) | Maximum Gas Flow Violation $\xi$ |
|----|------------------|------------------------|-----------------------|----------------------------------|
| M1 | 16750.84         | 4904.97                | 2776.36               | $1.68 \times 10^{-5}$            |
| M2 | 16750.84         | 4904.97                | 2776.36               | 0                                |
| M3 | 16935.74         | 4583.41                | 3107.57               | $8.56 \times 10^{-8}$            |
| M4 | 16750.85         | 4904.96                | 2776.37               | $7.61 \times 10^{-6}$            |

*kcf is short for kilo cubic feet.

TABLE III
AVERAGE COMPUTATIONAL EFFICIENCY OF 50 RANDOM SCENARIOS FOR THE STEADY STATE OPTIMAL GAS FLOW.

|                          | Method | Time (s) | Average Iteration |
|--------------------------|--------|----------|-------------------|
| 7-node natural gas system | M1 | 0.0782 | —— |
|                          | M2 | 0.0678 | —— |
|                          | M3 | 0.0623 | 2.00 |
|                          | M4 | 0.0423 | 1.06 |
| 20-node natural gas system | M1 | 0.5865 | —— |
|                          | M2 | 0.2164 | —— |
|                          | M3 | 0.7911 | 14.17 |
|                          | M4 | 0.0550 | 1.47 |

TABLE IV
ACCURACY AND EFFICIENCY OF 50 RANDOM SCENARIOS FOR THE QUASI-DYNAMIC OPTIMAL GAS FLOW OF A 20-NODE GAS SYSTEM.

|    | Time (s) | Average Iteration | Average Optimality Gap* | Average Violation $\xi$ |
|----|----------|-------------------|-------------------------|-------------------------|
| M1 | 29.88    | ——                | $4.93 \times 10^{-3}$%  | $1.30 \times 10^{-2}$   |
| M2 | 1.4935   | ——                | ——                      | 0                       |
| M3 | 3.29     | 7.92              | $1.02 \times 10^{-2}$%  | $4.93 \times 10^{-4}$   |
| M4 | 0.3326   | 1.78              | $4.89 \times 10^{-4}$%  | $3.28 \times 10^{-4}$   |

* The nonlinear solver (M2) is selected as the benchmark method for calculating the average optimality gap.

optimality gap. On the other hand, the data-driven CCP only requires 1.78 iterations before convergence on average, with 46% of the scenarios requiring only 1 iteration and 38% of the scenarios requiring 2. It therefore outperforms the cold-start CCP which requires 7.92 iterations on average. In terms of optimality, data-driven CCP is also preferable, with an average optimality gap of merely $4.89 \times 10^{-4}$%.

## V. CONCLUSION

We showed that data-driven approaches can boost the efficiency of convex optimization, especially SOCR, in the optimal gas flow problem. ANN can predict the nodal pressures with high precision, which will be set as a warm start for the iterative relaxations. The case studies, by comparing it with other model-driven methods, demonstrate its computational efficiency, feasibility and optimality.